\documentclass[12pt]{article}
\usepackage{mathrsfs,graphicx,amsmath,amssymb,amsfonts}
\begin{document}
\title{Some graft transformations and their applications on distance (signless) Laplacian spectra of graphs \thanks{This work is supported by NSFC (No. 11461071).}}
\author{Dandan Fan, Guoping Wang\footnote{Corresponding author. Email: xj.wgp@163.com.}\\
{\small School of Mathematical Sciences, Xinjiang Normal University,}\\
{\small Urumqi, Xinjiang 830054, P.R.China}}
\date{}
\maketitle {\bf Abstract.} Suppose that the vertex set of a connected graph $G$ is $V(G)=\{v_1,\cdots,v_n\}$.
Then we denote by $Tr_{G}(v_i)$ the sum of distances between $v_i$ and all other vertices of $G$.
Let $Tr(G)$ be the $n\times n$ diagonal matrix with its $(i,i)$-entry equal to $Tr_{G}(v_{i})$
and $D(G)$ be the distance matrix of $G$. Then $Q_{D}(G)=Tr(G)+D(G)$ and $L_{D}(G)=Tr(G)-D(G)$ are respectively the distance signless Laplacian matrix and distance Laplacian matrix of $G$.
The largest eigenvalues $\rho_{Q}(G)$ and $\rho_{L}(G)$ of $Q_D(G)$ and $L_D(G)$ are respectively called distance signless Laplacian spectral radius and distance Laplacian spectral radius of $G$.
In this paper we give some graft transformations and use them to characterize the tree $T$ such that $\rho_{Q}(T)$ and $\rho_{L}(T)$ attain the maximum
among all trees of order $n$ with given number of pendant vertices.

\section{Introduction}
Let $G$ be a simple connected graph with the vertex set $V(G)=\{v_1,\ldots,v_n\}$.
The {\it distance} $d_{G}(u,v)$ between the vertices $u$ and $v$ is length of a shortest path between $u$ and $v$ in $G$.
For $u\in V(G)$, the transmission $Tr_{G}(u)$ of $u$ is the sum of distances between $u$ and all other vertices of $G$.
Let $Tr(G)$ be the $n\times n$ diagnonal matrix with its $(i,i)$-entry equal to $Tr_{G}(v_{i})$ and $D(G)$ be the distance matrix of $G$.
Then $Q_{D}(G)=Tr(G)+D(G)$ and $L_{D}(G)=Tr(G)-D(G)$ are the distance signless Laplacian matrix and Laplacian matrix of $G$ respectively.
The largest eigenvalues $\rho_{Q}(G)$ and $\rho_{L}(G)$ of $Q_{D}(G)$ and $L_{D}(G)$ are distance signless Laplacian spectral radius and distance Laplacian spectral radius of $G$, respectively.

The distance spectral radius of a connected graph has been studied extensively.
S. Bose, M. Nath and S. Paul \cite{S. Bose} determined the unique graph with maximal distance
spectral radius in the class of graphs without a pendant vertex.
G. Yu, Y. Wu, Y. Zhang and J. Shu \cite{G. Yu-2} obtained respectively the extremal graph and unicyclic graph with the maximum and minimum distance spectral radius.
W. Ning, L. Ouyang and M. Lu \cite{W. Ning} characterized the graph with minimum distance spectral radius among trees with given number of pendant vertices.
G. Yu, H. Jia, H. Zhang and J. Shu \cite{G. Yu} characterized the unique graphs with the minimum and maximum distance spectral radius among trees and graphs with given number of pendant vertices.
For more about the distance spectra of graphs see the surveys \cite{M. Aouchiche-2, D. Stevanovic} as well as the references therein.

M. Aouchiche and P. Hansen introduced in \cite{M. Aouchiche} the distance Laplacian and distance signless Laplacian spectra of graphs,
and proved in \cite{M. Aouchiche-1} that the star $S_n$ attains minimum distance
Laplacian spectral radius among all trees of order $n$.
R. Xing and B. Zhou \cite{R. Xing} gave the unique graph with minimum distance and
distance signless Laplacian spectral radii among bicyclic graphs.
R. Xing, B. Zhou and J. Li \cite{R. Xing-2} determined the graphs with minimum
distance signless Laplacian spectral radius among the trees, unicyclic graphs,
bipartite graphs, the connected graphs with fixed pendant vertices and fixed connectivity, respectively.
A. Niu, D. Fan and G. Wang \cite{A. Niu} determined the extremal graph with the minimum distance Laplacian spectral radius among bipartite graphs with the fixed connectivity and matching numbers.
H. Lin and B. Zhou \cite{H. Lin-2} characterized the unique graph with the minimum distance Laplacian spectral radius among graphs with fixed number of pendant
vertices and edge connectivity.

In this paper we give some graft transformations and use them to characterize the tree $T$ such that $\rho_{Q}(T)$ and $\rho_{L}(T)$ attain the maximum
among all trees of order $n$ with given number of pendant vertices.

\section{Some graft transformations and their applications on distance Laplacian spectra}
If $x=(x_{1},x_{2},\ldots,x_{n})^T$ then it can be considered as a function defined on the vertex set $V(G)=\{v_{1},v_{2},\ldots,v_{n}\}$ of a graph
$G$ which maps vertex $v_{i}$ to $x_{v_i}$, and so
$$x^{T}L_{D}(G)x=\sum\limits_{\{u,v\}\subseteq V(G)}d_G(u,v)(x(u)-x(v))^{2}$$
which shows that $L_{D}(G)$ is positive semidefinite.

Suppose that $x$ is an eigenvector of $L_{D}(G)$ corresponding to the eigenvalue $\mu$.
Then for each $v\in V(G)$, $\mu x(v)=\sum\limits_{u\in V(G)}d_G(u,v)(x(v)-x(u))$.
Clearly, $\overrightarrow{1}=(1,1,\ldots,1)^{T}$ is an eigenvector of $L_{D}(G)$ corresponding to the eigenvalue $0$.

The union $G_{1}\cup G_{2}$ of graph $G_{1}$ and $G_{2}$ is the graph with the vertex set $V(G_{1})\cup V(G_{2})$ and the edge set $E(G_{1})\cup E(G_{2})$.\vskip 2mm

Now we give some graft transformations on distance Laplacian spectra of graphs. \vskip 2mm

{\bf Lemma 2.1.} {\it Let $G=G_{1}\cup G_{2}\cup G_{3}$, where $V(G_{i})\cap V(G_{j})=v_{0}$ for $1\leq i\neq j\leq 3$ and $V(G_{i})\geq 2$ for $i=1,2,3$.
Suppose that $u\in V(G_{2})\backslash \{v_{0}\}$.
Let $\widetilde{G}$ be the graph obtained from $G$ by moving $G_3$ from $v_0$ to $u$.
Let $x=(x_{v_0},x_{v_1},\ldots,x_{v_{n-1}})^T$ be a unit eigenvector corresponding to $\rho_{L}(G)$, where $x_{v_i}=x(v_i)$.
Then we have the following
\begin{enumerate}
\setlength{\parskip}{0ex}
\item[\rm (i)]
$\rho_{L}(\widetilde{G})\geq \rho_{L}(G)$ if $\sum\limits_{v_{i}\in V(G_{3})\setminus \{v_{0}\}}\sum\limits_{v_{j}\in V(G_{1})}(x_{v_i}-x_{v_j})^{2}\geq \sum\limits_{v_{i}\in V(G_{3})\setminus \{v_{0}\}}\sum\limits_{v_{j}\in V(G_{2})}(x_{v_i}-x_{v_j})^{2}$.
\item[\rm (ii)]
$\rho_{L}(\widetilde{G})>\rho_{L}(G)$ if $\sum\limits_{v_{i}\in V(G_{3})\setminus \{v_{0}\}}\sum\limits_{v_{j}\in V(G_{1})}(x_{v_i}-x_{v_j})^{2}> \sum\limits_{v_{i}\in V(G_{3})\setminus \{v_{0}\}}\sum\limits_{v_{j}\in V(G_{2})}(x_{v_i}-x_{v_j})^{2}$.
\end{enumerate}}

{\bf Proof.} Let $d_0=d_G(u,v_0)$. Then $d_0=d_{\widetilde{G}}(u,v_0)$ and for $v_i\in V(G_{3})\backslash \{v_0\}, v_{j}\in V(G_{1})$, $d_{\widetilde{G}}(v_{i},v_{j})=d_{G}(v_{i},v_{j})+d_{0}$;
and for $v_{i}\in V(G_{3})\backslash \{v_{0}\}, v_{j}\in V(G_{2})$,
$d_{\widetilde{G}}(v_{i},v_{j})\leq d_{G}(v_{i},v_{j})+ d_{0}$.

Thus, by Rayleigh's inequalities, we have

$\rho_{L}(\widetilde{G})-\rho_{L}(G)\geq x^{T}(L_{D}(\widetilde{G})-L_D(G))x$

\qquad\qquad\qquad\qquad\:\:\:\,\,\:$\geq d_{0}\sum\limits_{v_{i}\in V(G_{3})\setminus \{v_{0}\}}(\sum\limits_{v_{j}\in V(G_{1})}(x_{v_i}-x_{v_j})^{2}-\sum\limits_{v_{j}\in V(G_{2})}(x_{v_i}-x_{v_j})^{2})$\\
from which we can see that both (i) and (ii) are true. $\Box$\vskip 2mm

Denote by $N_G(v)$ and by $d_{G}(v)$ the neighbours and the degree of the vertex $v$ in graph $G$, respectively.
Denote by $Tr_{max}(G)$ the maximum transmission of vertices of $G$. In \cite{M. Nath} M. Nath and S. Paul pointed out that $\rho_{L}(G)>Tr_{max}(G)+1$.\vskip 2mm

{\bf Lemma 2.2.} {\it Suppose that $P=v_1v_2\ldots v_{\ell}$ is a path of order $\ell$, and $S_j$ be a star with the center $u_j$ $(j=1, i, \ell; 2\leq i\leq \ell-1)$.
Let $H$ be a graph obtained from $P$ by attaching the center $u_j$ of $S_j$ with $v_j$ for $j=1, i, \ell$ as shown in Fig. 1,
and $H_t$ be the graph obtained from $H$ by moving $S_i$ from $v_i$ to $v_t$ $(t=1~\mbox{or}~\ell)$. Then $\mbox{max}\{\rho_{L}(H_1), \rho_{L}(H_{\ell})\}>\rho_{L}(H)$.}
\begin{center}
\begin{figure}[htbp]
\centering
\includegraphics[height=25mm]{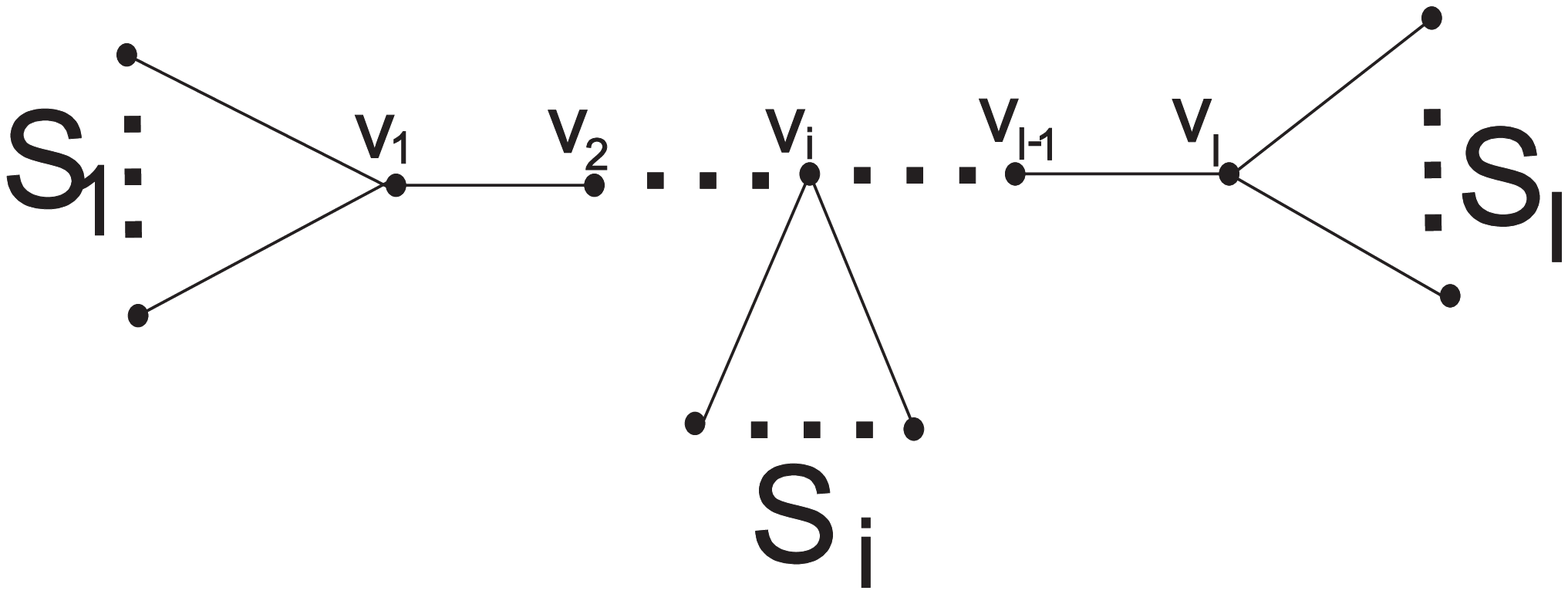}
\\ Fig. 1. ~~$H$
\end{figure}
\end{center}

{\bf Proof.} Let $x=(x_{v_{1}},x_{v_{2}},\cdots ,x_{v_{n}})^T$ be a unit eigenvector of
$H$ corresponding to $\rho_{L}(H)$ orthogonal to $\overrightarrow{1}$.
Set $\widetilde{S}_{1}=V(S_{1})\cup\{v_{2},v_{3},\ldots,v_{i-1}\}$, $\widetilde{S_{i}}=V(S_{i})\backslash \{v_{i}\}$
and $\widetilde{S}_{\ell}=V(S_{\ell})\cup\{v_{i},v_{i+1},\ldots,v_{\ell-1}\}$.

We first take $i=2$.
If $\sum\limits_{v_{m}\in \widetilde{S}_{2}}\sum\limits_{v_{j}\in \widetilde{S}_{1}}(x_{v_m}-x_{v_j})^{2}> \sum\limits_{v_{m}\in \widetilde{S}_{2}}\sum\limits_{v_{j}\in \widetilde{S}_{\ell}}(x_{v_m}-x_{v_j})^{2}$,
then by Lemma 2.1 (ii), we have $\rho_{L}(H_{\ell})>\rho_{L}(H)$.
So we assume $\sum\limits_{v_{m}\in \widetilde{S}_{2}}\sum\limits_{v_{j}\in \widetilde{S}_{\ell}}(x_{v_m}-x_{v_j})^{2}\geq \sum\limits_{v_{m}\in \widetilde{S}_{2}}\sum\limits_{v_{j}\in \widetilde{S}_{1}}(x_{v_m}-x_{v_j})^{2}$.
Then by Lemma 2.1 (i), we have $\rho_{L}(H_1)\geq\rho_{L}(H)$.

Now we suppose for a contradiction that $\rho_{L}(H)=\rho_{L}(H_1)$. Thus, $x$ is also a unit eigenvector of $H_1$ corresponding to $\rho_{L}(H_1)$.
Then, for any $v_{m}\in\widetilde{S}_{1}$, we have
$$\left\{
  \begin{array}{ll}
 \rho_{L}(H)x_{v_{m}}=\sum\limits_{v_{j}\in \widetilde{S}_{1}}d_{H}(v_{m},v_{j})(x_{v_{m}}-x_{v_{j}})+\sum\limits_{v_{j}\in \widetilde{S}_{2}}d_{H}(v_{m},v_{j})(x_{v_{m}}-x_{v_{j}})& \\
 \qquad\qquad\qquad+\sum\limits_{v_{j}\in\widetilde{S}_{\ell}}d_{H}(v_{m},v_{j})(x_{v_{m}}-x_{v_{j}}),& \\
 \rho_{L}(H_1)x_{v_{m}}=\sum\limits_{v_{j}\in \widetilde{S}_{1}}d_{H_1}(v_{m},v_{j})(x_{v_{m}}-x_{v_{j}})+\sum\limits_{v_{j}\in \widetilde{S}_{2}}d_{H_1}(v_{m},v_{j})(x_{v_{m}}-x_{v_{j}})& \\
 \qquad\qquad\qquad+\sum\limits_{v_{j}\in \widetilde{S}_{\ell}}d_{H_1}(v_{m},v_{j})(x_{v_{m}}-x_{v_{j}}).&  \hbox{}
  \end{array}
\right.$$
Since $d_{H}(v_{m},v_{j})=d_{H_1}(v_{m},v_{j})$ for $v_{j}\in \widetilde{S}_{1}$, $d_{H}(v_{m},v_{j})=d_{H_1}(v_{m},v_{j})$ for $v_{j}\in\widetilde{S}_{\ell}$ and
$d_{H}(v_{m},v_{j})=d_{H_1}(v_{m},v_{j})+(i-1)$ for $v_{j}\in\widetilde{S}_{2}$, we can obtain from the above equations that
$$(\rho_{L}(H)-\rho_{L}(H_1))x_{v_{m}}=(i-1)\sum\limits_{v_{j}\in \widetilde{S}_{2}}(x_{v_{m}}-x_{v_{j}})=0.$$

Let $x_{\widetilde{S}_{1}}=\sum\limits_{v\in \widetilde{S}_{1}}x_{v}$, $x_{\widetilde{S}_{2}}=\sum\limits_{v\in \widetilde{S}_{2}}x_v$ and $x_{\widetilde{S}_{\ell}}=\sum\limits_{v\in\widetilde{S}_{\ell}}x_v$.
Then we have $x_{v_{m}}=\frac{x_{\widetilde{S}_{2}}}{|\widetilde{S}_{2}|}$.

Similarly, we have $x_{v_{k}}=\frac{x_{\widetilde{S}_{2}}}{|\widetilde{S}_{2}|}$ for any $v_{k}\in \widetilde{S}_{\ell}$.

We can see that on the one hand, $\sum\limits_{j=1}^{n}x_{v_{j}}=x_{\widetilde{S}_{1}}+x_{\widetilde{S}_{2}}+x_{\widetilde{S}_{\ell}}=\frac{|\widetilde{S}_{1}|x_{\widetilde{S}_{2}}}{|\widetilde{S}_{2}|}+\frac{|\widetilde{S}_{\ell}|x_{\widetilde{S}_{2}}}{|\widetilde{S}_{2}|}
+x_{\widetilde{S}_{2}}=(\frac{|\widetilde{S}_{1}|+|\widetilde{S}_{2}|+|\widetilde{S}_{\ell}|}{|\widetilde{S}_{2}|})
x_{\widetilde{S}_{2}}$, and on the other hand, $\sum\limits_{j=1}^{n}x_{v_{j}}=0$, and so
$x_{\widetilde{S}_{2}}=0$.

If $|\widetilde{S}_{\ell}|\neq|\widetilde{S}_{1}|$ then we have $x_{v_{p}}=\frac{x_{\widetilde{S}_{\ell}}-x_{\widetilde{S}_{1}}}{|\widetilde{S}_{\ell}|-|\widetilde{S}_{1}|}$ for any $v_{p}\in \widetilde{S}_{2}$,
and so $x(v)=0$ for any $v\in \widetilde{S}_{2}$,
which implies that $x=\{0,0,\ldots,0\}^{T}$.
This contradiction shows that $\rho_{L}(H_1)>\rho_{L}(H)$.

If $|\widetilde{S}_{\ell}|=|\widetilde{S}_{1}|$ then we easily observe that $Tr_{max}(H_1)=Tr_{H_1}(v_k)$,
where $v_k$ is a pendent vertex in $S_\ell$.
For any $v_{q}, v_{p}\in \widetilde{S}_{2}$, we have
$$\left\{
  \begin{array}{ll}
 \rho_{L}(H_1)x_{v_{q}}=\sum\limits_{v_{j}\in \widetilde{S}_{1}}d_{H_1}(v_{q},v_{j})(x_{v_{q}}-x_{v_{j}})+\sum\limits_{v_{j}\in \widetilde{S}_{2}}d_{H_1}(v_{q},v_{j})(x_{v_{q}}-x_{v_{j}})& \\
 \qquad\qquad\qquad+\sum\limits_{v_{j}\in\widetilde{S}_{\ell}}d_{H_1}(v_{q},v_{j})(x_{v_{q}}-x_{v_{j}}),& \\
 \rho_{L}(H_1)x_{v_{p}}=\sum\limits_{v_{j}\in \widetilde{S}_{1}}d_{H_1}(v_{p},v_{j})(x_{v_{p}}-x_{v_{j}})+\sum\limits_{v_{j}\in \widetilde{S}_{2}}d_{H_1}(v_{p},v_{j})(x_{v_{p}}-x_{v_{j}})& \\
 \qquad\qquad\qquad+\sum\limits_{v_{j}\in \widetilde{S}_{\ell}}d_{H_1}(v_{p},v_{j})(x_{v_{p}}-x_{v_{j}}).&  \hbox{}
  \end{array}
\right.$$
from which we obtain that $ (\rho_{L}(H_1)-(Tr_{H_1}(v_p)+2))(x_{v_{q}}-x_{v_{p}})=0$.
We know that $\rho_{L}(H_1)> Tr_{H_1}(v_k)+1$.
By a simple computation we can know that $Tr_{H_1}(v_k)>Tr_{H_1}(v_p)+2$, and so
$x_{v_{q}}=x_{v_{p}}$. But $x_{\widetilde{S}_{2}}=0$, and so $x(v)=0$ for any $v\in \widetilde{S}_{2}$,
which implies that $x=\{0,0,\ldots,0\}^{T}$.
This contradiction shows that $\rho_{L}(H_1)>\rho_{L}(H)$.

For the case $i=\ell-1$, we can similarly prove that $\rho_{L}(H_1)>\rho_{L}(H)$ and $\rho_{L}(H_{\ell})>\rho_{L}(H)$.
Next we suppose that $3\leq i\leq \ell-2$.
If $\sum\limits_{v_{m}\in \widetilde{S}_{i}}\sum\limits_{v_{j}\in \widetilde{S}_{\ell}}(x_{v_m}-x_{v_j})^{2}\geq \sum\limits_{v_{m}\in \widetilde{S}_{i}}\sum\limits_{v_{j}\in \widetilde{S}_{1}}(x_{v_m}-x_{v_j})^{2}$,
then let $H'=H-\sum\limits_{v_{j}\in \widetilde{S}_i}v_{i}v_{j}+\sum\limits_{v_{j}\in \widetilde{S}_i}v_{2}v_j$,
and otherwise, let $H'=H-\sum\limits_{v_{j}\in \widetilde{S}_i}v_{i}v_{j}+\sum\limits_{v_{j}\in \widetilde{S}_i}v_{\ell-1}v_j$.
By Lemma 2.1, $\rho_{L}(H')\geq\rho_{L}(H)$.
The former argument has shown that $\mbox{max}\{\rho_{L}(H_1), \rho_{L}(H_{\ell})\}>\rho_{L}(H')$,
and so $\mbox{max}\{\rho_{L}(H_1), \rho_{L}(H_{\ell})\}>\rho_{L}(H)$. $\Box$\vskip 2mm

A path $u_{1}\ldots u_{r}$ ($r\geq2$) in a graph $G$ is a {\it pendant path} at $u_1$ if $d_{G}(u_1)\geq3$, $d_{G}(u_{i})=2$ ($2\leq i\leq r-1$) and $d_{G}(u_{r})=1$.
If graphs $X$ and $Y$ are isomorphic then we write $X\cong Y$, and otherwise $X\not\cong Y$.\vskip 2mm

{\bf Lemma 2.3.} {\it Let $G_{p,q}$ be the graph obtained from a connected graph $G$ of order at least two
by attaching two pendant paths $P_1$ and $P_2$ of length $p$ and $q$ respectively at a vertex $u$ of $G$.
If $q\geq p\geq1$, then $\rho_{L}(G_{p-1,q+1})>\rho_{L}(G_{p,q})$.} \vskip 2mm

{\bf Proof.} Suppose that $P_{1}=uu_{1}u_{2}\ldots u_{p}$ and $P_{2}=uv_{1}v_{2}\ldots v_{q}$.
Let $S=V(G)\backslash\{u\}$ and $P'_{i}=V(P_{i})\backslash\{u\}$ ($i=1,2$).
Let $x$  be a unit eigenvector of $G_{p,q}$ corresponding to $\rho_{L}(G_{p,q})$ orthogonal to $\overrightarrow{1}$.

If $q>p$, we suppose on the contrary that $\sum\limits_{v_{i}\in P'_{2}}\sum\limits_{v_{k}\in S}(x_{v_i}-x_{v_k})^{2}< \sum\limits_{v_{i}\in P'_{1}}\sum\limits_{v_{k}\in S}(x_{v_i}-x_{v_k})^{2}$.

Since $G_{p,q}\cong G_{q,p}$, it follows from Rayleigh's inequalities that

$0=\rho_{L}(G_{q,p})-\rho_{L}(G_{p,q})\,\:\geq x^{T}(L_D(G_{q,p})-L_D(G_{p,q}))x$

\qquad\qquad\qquad\qquad\qquad\quad$\geq(q-p)(\sum\limits_{v_{i}\in V(P_{1})}\sum\limits_{v_{k}\in S}(x_{v_i}-x_{v_k})^{2}- \sum\limits_{v_{i}\in P'_{2}}\sum\limits_{v_{k}\in S}(x_{v_i}-x_{v_k})^{2})$

\qquad\qquad\qquad\qquad\qquad\quad$\geq(q-p)(\sum\limits_{v_{i}\in P'_{1}}\sum\limits_{v_{k}\in S}(x_{v_i}-x_{v_k})^{2}- \sum\limits_{v_{i}\in P'_{2}}\sum\limits_{v_{k}\in S}(x_{v_i}-x_{v_k})^{2})>0$.

This contradiction shows that $\sum\limits_{v_{i}\in P'_{2}}\sum\limits_{v_{k}\in S}(x_{v_i}-x_{v_k})^{2}\geq \sum\limits_{v_{i}\in P'_{1}}\sum\limits_{v_{k}\in S}(x_{v_i}-x_{v_k})^{2}$.

Now we suppose that $p=q$.
Since $P_1$ and $P_2$ have the same order, we can exchange the labels of $P_1$ and $P_2$,
and so we can assume that $\sum\limits_{v_{i}\in P'_{2}}\sum\limits_{v_{k}\in S}(x_{v_i}-x_{v_k})^{2}\geq \sum\limits_{v_{i}\in P'_{1}}\sum\limits_{v_{k}\in S}(x_{v_i}-x_{v_k})^{2}$.

Note that

$\rho_{L}(G_{p-1,q+1})-\rho_{L}(G_{p,q})\:\,\geq x^{T}(L_D(G_{p-1,q+1})-L_{D}(G_{p,q}))x$

\qquad\qquad\qquad\qquad\qquad\quad$\geq\sum\limits_{v_{i}\in V(P_{2})}\sum\limits_{v_{k}\in S}(x_{v_i}-x_{v_k})^{2}- \sum\limits_{v_{i}\in P'_{1}}\sum\limits_{v_{k}\in S}(x_{v_i}-x_{v_k})^{2}$

\qquad\qquad\qquad\qquad\qquad\quad$\geq\sum\limits_{v_{i}\in P'_{2}}\sum\limits_{v_{k}\in S}(x_{v_i}-x_{v_k})^{2}- \sum\limits_{v_{i}\in P'_{1}}\sum\limits_{v_{k}\in S}(x_{v_i}-x_{v_k})^{2}$.\\

Therefore, $\rho_{L}(G_{p-1,q+1})\geq\rho_{L}(G_{p,q})$.

Now we suppose for a contradiction that $\rho_{L}(G_{p-1,q+1})=\rho_{L}(G_{p,q})$. Thus, $x$ is also a unit eigenvector of $G_{p-1,q+1}$ corresponding to $\rho_{L}(G_{p-1,q+1})$.
Then, for any $v_{m}\in P'_{1}$, we have
$$\left\{
  \begin{array}{ll}
 \rho_{L}(G_{p,q})x_{v_{m}}=\sum\limits_{v_{i}\in P'_{1}}d_{G_{p,q}}(v_{m},v_{i})(x_{v_{m}}-x_{v_{i}})+\sum\limits_{v_{i}\in V(P_{2})}d_{G_{p,q}}(v_{m},v_{i})(x_{v_{m}}-x_{v_{i}})& \\
 \qquad\qquad\qquad+\sum\limits_{v_{i}\in S}d_{G_{p,q}}(v_{m},v_{i})(x_{v_{m}}-x_{v_{i}}),& \\
 \rho_{L}(G_{p-1,q+1})x_{v_{m}}=\sum\limits_{v_{i}\in P'_{1}}d_{G_{p-1,q+1}}(v_{m},v_{i})(x_{v_{m}}-x_{v_{i}})+\sum\limits_{v_{i}\in V(P_{2})}d_{G_{p-1,q+1}}(v_{m},v_{i})(x_{v_{m}}-x_{v_{i}})& \\
 \qquad\qquad\qquad+\sum\limits_{v_{i}\in S}d_{G_{p-1,q+1}}(v_{m},v_{i})(x_{v_{m}}-x_{v_{i}}).&  \hbox{}
  \end{array}
\right.$$
Since $d_{G_{p-1,q+1}}(v_{m},v_{i})=d_{G_{p,q}}(v_{m},v_{i})$ for $v_{i}\in P'_{1}$, $d_{G_{p-1,q+1}}(v_{m},v_{i})=d_{G_{p,q}}(v_{m},v_{i})$ for $v_{i}\in V(P_{2})$ and
$d_{G_{p-1,q+1}}(v_{m},v_{i})=d_{G_{p,q}}(v_{m},v_{i})-1$ for $v_{i}\in S$, we can obtain that

$(\rho_{L}(G_{p-1,q+1})-\rho_{L}(G_{p,q}))x_{v_{m}}=- \sum\limits_{v_{i}\in S}(x_{v_{m}}-x_{v_{i}})=0.$

Let $x_{P'_{1}}=\sum\limits_{v\in P'_{1}}x_{v}$, $x_{P_{2}}=\sum\limits_{v\in V(P_{2})}x_v$ and $x_{S}=\sum\limits_{v\in S}x_v$.
Note that $|S|=n-p-q-1$.
Thus, from the above equation, we have $x_{v_{m}}=\frac{x_{S}}{n-p-q-1}$.

Similarly, we have $x_{v_{k}}=\frac{x_{S}}{n-p-q-1}$ for $v_{k}\in V(P_{2})$
and $x_{v_{p}}=\frac{x_{P'_{1}}-x_{P_{2}}}{p-q-1}$ for $v_{p}\in S$.

We can see that on the one hand, $\sum\limits_{v\in{G_{p,q}}}x_{v}=x_{P'_{1}}+x_{P_{2}}+x_{S}=\frac{nx_{S}}{n-p-q-1}$, and on the other hand, $\sum\limits_{v\in{G_{p,q}}}x_{v}=0$, and so
$x_{S}=0$, which implies that $x=\{0,0,\ldots,0\}^{T}$.
This contradiction shows that $\rho_{L}(G_{p-1,q+1})>\rho_{L}(G_{p,q})$. $\Box$\vskip 2mm

As one application of above transformations we next determine the tree and graph that attain the maximum distance Laplacian spectral radius
among all trees and graphs of order $n$ with given number of pendant vertices, respectively.\vskip 2mm

Denote by $T(n,k; t_1,t_2)$ the graph obtained from a path $P_{\ell}$ of order $\ell$ ($\ell\geq 2$) by attaching $t_1$ edges at an end of $P_{\ell}$ and $t_2$ edges at other end
of $P_{\ell}$, where $t_{1}+t_{2}=k$ and $\ell+k=n$. Clearly, $P_n\cong T(n,2; 1,1)$.
Let $\mathfrak{D}(n,k)=\{T(n,k; t_1,t_2)|~ t_{1}\geq 1, t_{2}\geq 1~ \mbox{and}~ t_{1}+t_{2}=k\geq 2\}.$
Let $\mathfrak{T}(n,k)=\{T|~T~ \mbox{is a tree of order}~n~ \mbox{with}~ k ~\mbox{pendant vertices}, n-2 \geq k \geq 2\}.$\vskip 2mm

{\bf Theorem 2.4.} {\it The tree which attains the maximum distance Laplacian spectral radius in $\mathfrak{T}(n,k)$ belongs to $\mathfrak{D}(n,k)$.}

{\bf Proof.} Suppose that $T$ is a tree in $\mathfrak{T}(n,k)$, and let $n^3_T$ be the number of vertices which are of degree at least three in $T$.
If $n^3_T=0$ then $T\cong T(n,2; 1,1)$.
So we assume $n^3_T\geq 1$, and let $A=\{w_{1},w_{2},\ldots,w_{n^3_T}\}$ be the set of all vertices in $T$ which are of degree at least three.
Next we distinguish the following three cases to discuss.

{\bf Case 1.} $n^3_T=1$.

If $T\cong T(n,k; 1,k-1)$ then the proof is finished,
and otherwise $T$ must contain two pendant paths $T_1$ and $T_2$ of order at least three.
Suppose $T_{1}=w_{1}u_{1} \ldots u_{s}$ ($s\geq2$) and $T_{2}=w_{1}v_{1} \ldots v_{t}$ ($t\geq2$),
and let $T=T_{1}\cup T_{2}\cup T_{3}$, where $V(T_{i})\cap V(T_{j})=w_1$ for $1\leq i\neq j\leq 3$.
Without loss of generality we can assume that $s\geq t$.
Let $T'=T-v_{t-1}v_{t}+u_{s}v_{t}$. Then $T'\in\mathfrak{T}(n,k)$, and by Lemma 2.3, $\rho_{L}(T')>\rho_{L}(T)$.
Continuing the above process, we can finally see that $\rho_{L}(T(n,k; 1,k-1))>\rho_{L}(T)$.

{\bf Case 2.} $n^3_T=2$.

In this case $T$ contains exactly two vertices $w_1$ and $w_2$ of degrees at least three.
Suppose that $P=w_{1}v_{1}\ldots v_{t}w_{2}$ is the internal path in $T$.
If $T\in \mathfrak{D}(n,k)$ then the proof is finished,
and otherwise $T$ must contain a pendant path $T_2$ of order at least three.
Suppose that $T_{2}=w_{1}u_{1} \ldots u_{s}$ ($s\geq2$), and
let $T=T_{1}\cup T_{2}\cup T_{3}$, where $V(T_{i})\cap V(T_{j})=w_1$ for $1\leq i\neq j\leq 3$, $N_{T}(w_{1})\backslash\{u_{1},v_{1}\} \subset V(T_{1})$ and $w_{2}\in V(T_3)$, as in Fig. 2.
\begin{center}
\begin{figure}[htbp]
\centering
\includegraphics[height=20mm]{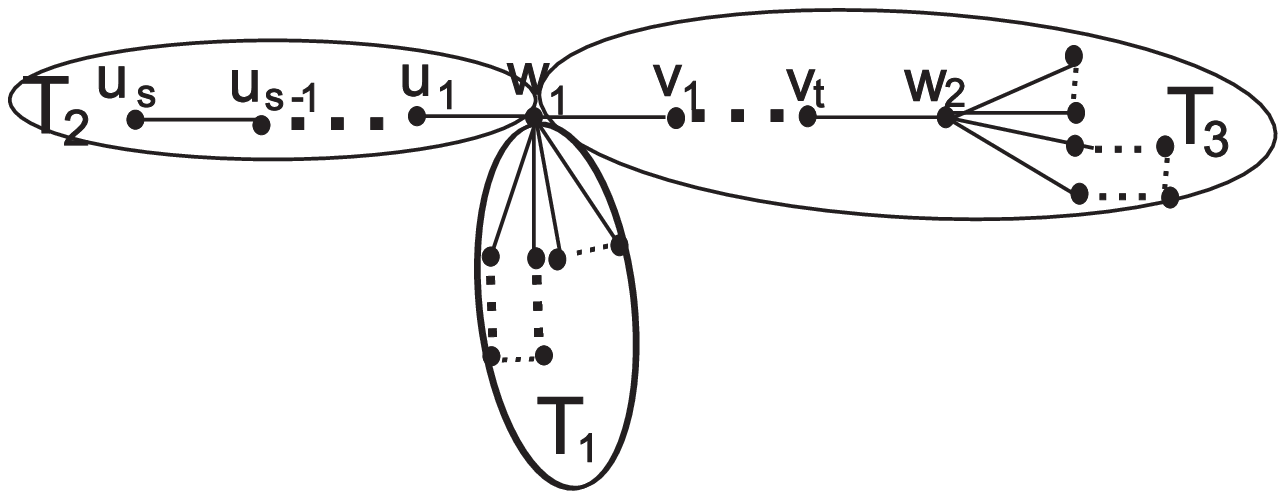}
\\ Fig. 2. ~~$T=T_{1}\cup T_{2}\cup T_{3}$
\end{figure}
\end{center}
If $\sum\limits_{v_{i}\in V(T_{1})}\sum\limits_{v_{j}\in V(T_{3})}(x_{v_i}-x_{v_j})^{2}\geq \sum\limits_{v_{i}\in V(T_{1})}\sum\limits_{v_{j}\in V(T_{2})}(x_{v_i}-x_{v_j})^{2}$,
then let $T'=T-\sum\limits_{v_{j}\in V(T_{1})}w_{1}v_{j}+\sum\limits_{v_{j}\in V(T_1)}u_{s-1}v_j$, and otherwise $T'=T-\sum\limits_{v_{j}\in V(T_{1})}w_{1}v_{j}+\sum\limits_{v_{j}\in V(T_1)}w_{2}v_j$.

Suppose that the former rises.
Then $T'\in\mathfrak{T}(n,k)$, but the number of pendant paths of length at least two in $T'$ is fewer than in $T$.
Continuing the above process, we finally get a resulting graph $T(n,k; t_1,t_2)$.
By the above argument, we know that $\rho_{L}(T(n,k; t_1,t_2))\geq \rho_{L}(T)$, with equality if and only if $T\cong T(n,k; t_1,t_2)$.

So we assume that the latter rises. Then $n^3_{T'}=1$.
If $T'\cong T(n,k; 1,k-1)$ then $T\cong H$, where $H$ is shown in Fig. 1.
By Lemma 2.2, $\rho_{L}(T(n,k;  1,k-1))>\rho_{L}(T)$.
If $T'\not\cong T(n,k; 1,k-1)$ then the Case 1 has shown that $\rho_{L}(T(n,k; 1,k-1))> \rho_{L}(T')$.
By Lemma 2.1, $\rho_{L}(T')\geq\rho_{L}(T)$, and so $\rho_{L}(T(n,k; 1,k-1))> \rho_{L}(T)$.

{\bf Case 3.} $n^3_T\geq 3$.

Without loss of generality we can assume that $d_{T}(w_1,w_2)=max\{d_{T}(w_p,w_q)\mid 1\leq p\neq q\leq n^3_T\}$.
Let $T=T_{1}\cup T_{2}\cup T_{3}$, where $V(T_{i})\cap V(T_{j})=w_{3}$ $(1\leq i\neq j\leq3)$, $w_{1}\in T_{1}$ and $w_{2}\in T_{2}$.

If $\sum\limits_{v_{i}\in V(T_{3})\setminus \{w_{3}\}}\sum\limits_{v_{j}\in V(T_{1})}(x_{v_i}-x_{v_j})^{2}\geq \sum\limits_{v_{i}\in V(T_{3})\setminus \{w_{3}\}}\sum\limits_{v_{j}\in V(T_{2})}(x_{v_i}-x_{v_j})^{2}$,
then let $T^1=T-\sum\limits_{v_{j}\in N_{T_{3}}(w_3)}w_{3}v_{j}+\sum\limits_{v_{j}\in N_{T_{3}}(w_3)}w_{2}v_{j}$,
and otherwise $T^1=T-\sum\limits_{v_{j}\in N_{T_{3}}(w_3)}w_{3}v_{j}+\sum\limits_{v_{j}\in N_{T_{3}}(w_3)}w_{1}v_{j}$.

It is obvious that $T^1\in \mathfrak{T}(n,k)$, but $n^3_{T^1}=n^3_{T}-1$.
After making some $r$ steps of the transformation, we can get a resulting graph $T^r\in \mathfrak{T}(n,k)$, but $n^3_{T^r}=2$.

If $T^r\in T(n,k; t_1,t_2)\in \mathfrak{D}(n,k)$, then $T^{r-1}\cong H$, where $H$ is shown in Fig. 1.
By Lemma 2.2, $\rho_{L}(T^r)>\rho_{L}(T^{r-1})$. Lemma 2.1 shows that $\rho_{L}(T^{r-1})\geq \rho_{L}(T)$, and so $\rho_{L}(T^r)>\rho_{L}(T)$.
If $T^r\not\in \mathfrak{D}(n,k)$ then the Case 2 has shown that there is some $T(n,k; t_1,t_2)\in \mathfrak{D}(n,k)$ such that $T(n,k; t_1,t_2)> \rho_{L}(T^r)$.
By Lemma 2.1, $\rho_{L}(T^r)\geq \rho_{L}(T)$, and so $T(n,k; t_1,t_2)> \rho_{L}(T)$.   $\Box$\vskip 2mm

{\bf Lemma 2.5} \cite{M. Aouchiche}. {\it If $u$ and $v$ are two non-adjacent vertices of graph $G$, then $\rho_{L} (G+uv)\leq \rho_{L}(G)$.}\vskip 2mm

Let $\Re(n,k)=\{G|~G~ \mbox{is a graph of order} ~n ~\mbox{with} ~k\geq 2~ \mbox{pendant vertices}\}.$\vskip 2mm

{\bf Theorem 2.6.} {\it The graph with the maximum distance Laplacian spectral radius in $\Re(n,k)$ lies to $\mathfrak{D}(n,k)$.}

{\bf Proof.} Suppose that $G\in \Re(n,k)$ and that $H$ is a spanning tree of $G$. Then by Lemma 2.5, $\rho_{L}(H)\geq\rho_{L}(G)$.
Let $\mathfrak{K}(H)$ denote the number of pendant vertices in the tree $H$.
Then $\mathfrak{K}(H)\geq k$.
Theorem 2.4 has told us that there is a tree $T(n,\mathfrak{K}(H); n_1,n_2)\in \mathfrak{D}(n,\mathfrak{K}(H))$
such that $\rho_{L}(T(n,\mathfrak{K}(H); n_1,n_2))\geq \rho_{L}(H)$.
If $\mathfrak{K}(H)=k$ then the proof is finished,
and otherwise we can construct a tree $T\in \mathfrak{T}(n,k)$ obtained from $H$
by deleting $\mathfrak{K}(H)-k$ pendant edges and then attaching them with the remaining pendant edges. By Lemma 2.3, $\rho_{L}(T)> \rho_{L}(T(n,\mathfrak{K}(H); n_1,n_2))$.
It follows from Theorem 2.4 that there is a tree $T(n,k; t_1,t_2)$ such that $\rho_{L}(T(n,k; t_1,t_2))> \rho_{L}(T)$ and so $\rho_{L}(T(n,k; t_1,t_2))>\rho_{L}(G)$.
$\Box$\vskip 2mm

As one consequence of the above argument, we directly have\vskip 2mm

{\bf Theorem 2.7.} {\it If $G\in \Re(n,3)$ attains the maximum distance Laplacian spectral radius then $G\cong T(n,3; 1,2)$.}

\section{Some graft transformations and their applications on distance signless Laplacian spectra}
Suppose that $G$ is a connected graph with $V(G)=\{v_{1},v_{2},\ldots,v_{n}\}$.
Let $x=(x_{v_1},x_{v_2},\ldots,x_{v_n})^T$, where $v_{i}$ corresponds to $x_{v_i}$ for $i=1,2,\ldots, n$.
Then $$x^{T}Q_{D}(G)x=\sum\limits_{\{u,v\}\subseteq V(G)}d_G(u,v)(x(u)+x(v))^{2}$$
which shows that $Q_{D}(G)$ is positive definite for $n\geq 3$.\vskip 2mm

Suppose that $x$ is an eigenvector of $Q_D(G)$ corresponding to the eigenvalue $\lambda$. Then for each $v\in V(G)$,
$\lambda x(v)=\sum\limits_{u\in V(G)}d_G(u,v)(x(u)+x(v)).$\vskip 2mm

If $x>0$ satisfies $Q_{D}(G)x=\rho_{Q}(G)x$ then we call $x$ the Perron vector of $Q_{D}(G)$ corresponding to $\rho_{Q}(G)$.\vskip 2mm

Now we give some graft transformations on distance signless Laplacian spectra of graphs.\vskip 2mm

{\bf Lemma 3.1.} {\it Let $G=G_{1}\cup G_{2}\cup G_{3}$, where $V(G_{i})\cap V(G_{j})=v_{0}$ for $1\leq i\neq j\leq 3$ and $V(G_{i})\geq 2$ for $i=1,2,3$.
Suppose that $u\in V(G_{2})\backslash \{v_{0}\}$, and
let $\widetilde{G}$ be the graph obtained from $G$ by moving $G_3$ from $v_0$ to $u$.
Let $x=(x_{v_{1}},x_{v_{2}},\cdots ,x_{v_{n}})^T$ be a Perron vector corresponding to $\rho_{Q}(G)$, where $x_{v_i}=x(v_i)$.
If $\sum\limits_{v_{i}\in V(G_{3})\setminus \{v_{0}\}}\sum\limits_{v_{j}\in V(G_{1})}(x_{v_i}+x_{v_j})^{2}\geq \sum\limits_{v_{i}\in V(G_{3})\setminus \{v_{0}\}}\sum\limits_{v_{j}\in V(G_{2})}(x_{v_i}+x_{v_j})^{2}$,
then $\rho_{Q}(\widetilde{G})>\rho_{Q}(G)$.}

{\bf Proof.} It can be easily observed that $d_G(u,v_0)=d_{\widetilde{G}}(u,v_0)$ and that for $v_i\in V(G_{3})\backslash \{v_0\}$, $v_{j}\in V(G_{1})$, $d_{\widetilde{G}}(v_{i},v_{j})=d_{G}(v_{i},v_{j})+d_G(u,v_0)$;
and for $v_{i}\in V(G_{3})\backslash \{v_{0}\}$, $v_{j}\in V(G_{2})$,
$d_{\widetilde{G}}(v_{i},v_{j})\leq d_{G}(v_{i},v_{j})+ d_G(u,v_0)$.

Thus, by Rayleigh's inequalities we have

$\rho_{Q}(\widetilde{G})-\rho_{Q}(G)\geq x^{T}(Q_{D}(\widetilde{G})-Q_D(G))x$

\qquad\qquad\quad\quad\qquad\quad\,\,\:\,\,$\geq d_G(u,v_0)\sum\limits_{v_{i}\in V(G_{3})\setminus \{v_{0}\}}(\sum\limits_{v_{j}\in V(G_{1})}(x_{v_i}+x_{v_j})^{2}-\sum\limits_{v_{j}\in V(G_{2})}(x_{v_i}+x_{v_j})^{2})$

\qquad\qquad\quad\quad\qquad\quad\,\,\:\,\,$\geq 0$.

Suppose on the contrary that $\rho_{Q}(\widetilde{G})=\rho_{Q}(G)$. Then $x^{T}Q_D(G)x=\rho_{Q}(\widetilde{G})$.

Therefore, $(\rho_{Q}(\widetilde{G})-\rho_{Q}(G))x_{v_0}=\sum\limits_{v_{k}\in V(G_{3})\backslash \{v_{0}\}}d_G(u,v_0)(x_{v_0}+x_{v_k})>0$.
This contradiction shows that $\rho_{Q}(\widetilde{G})>\rho_{Q}(G)$. $\Box$\vskip 2mm

{\bf Lemma 3.2} \cite{H. Liu}. {\it Let $G_{p,q}$ be the connected graph obtained from a graph $G$ of order at least two by attaching two pendant paths of length $p$ and $q$ at a vertex $u$ of $G$.
If $q\geq p\geq 1$, then $\rho_{Q}(G_{p-1,q+1})>\rho_{Q}(G_{p,q})$.} \vskip 2mm

As one application of above transformations we determine the tree and graph that attain the maximum distance signless Laplacian spectral radius in $\mathfrak{T}(n,k)$ and $\Re(n,k)$, respectively.\vskip 2mm

{\bf Theorem 3.3.} {\it The tree which attains the maximum distance signless Laplacian spectral radius in $\mathfrak{T}(n,k)$ belongs to $\mathfrak{D}(n,k)$.}

{\bf Proof.} Suppose that $T_0$ is a tree with the maximum distance signless Laplacian spectral radius in $\mathfrak{T}(n,k)$.
It is clear that if $k=2$ then $T_{0}\cong T(n,2; 1,1)$ and
that if $k=3$ then $T_{0}\cong T(n,3; 1,2)$. So we assume $k\geq 4$.
Let $n^3_{T_{0}}$ be the number of vertices which are of degree at least three in $T_0$.
Then $n^3_{T_{0}}\geq 1$.
Let $A=\{w_{1},w_{2},\ldots,w_{n^3_{T_{0}}}\}$ be the set of all vertices in $T_{0}$ which are of degree at least three.

Suppose that $n^3_{T_{0}}=1$. Then $d_{T_0}(w_1)\geq3$.
If $T_{0}\not\cong T(n,k; 1,k-1)$ then there must be two pendant paths of length at least two at $w_1$,
say, $P_{1}=w_{1}u_{1} \ldots u_{s}$ and $P_{2}=w_{1}v_{1} \ldots v_{t}$ ($s\geq t\geq2$).
Let $T'_{0}=T_{0}-v_{t-1}v_{t}+u_{s}v_{t}$. Then $T'_{0}\in\mathfrak{T}(n,k)$, and by Lemma 3.2, we have $\rho_{Q}(T'_{0})>\rho_{Q}(T_{0})$.
This contradicts the maximality of $T_{0}$, and so $T_{0}\cong T(n,k; 1,k-1)$.

Next we assume $n^3_{T_{0}}\geq 2$ and verify the following two claims.

{\bf Claim 1.} $n^3_{T_{0}}= 2$.

Suppose for a contradiction that $n^3_{T_{0}}\geq 3$. Then $w_{3}\in A$.
Let $d_{T_{0}}(w_1,w_2)=max\{d_{T_{0}}(w_p,w_q)\mid 1\leq p\neq q\leq n^3_{T_{0}}\}$.
Then we can let $T_{0}=T_{1}\cup T_{2}\cup T_{3}$,
where $V(T_{i})\cap V(T_{j})=w_{3}$ $(1\leq i\neq j\leq3)$ and $w_{1}\in T_{1}$, $w_{2}\in T_{2}$.

If $\sum\limits_{v_{i}\in V(T_{3})\setminus \{w_{3}\}}\sum\limits_{v_{j}\in V(T_{1})}(x_{v_i}-x_{v_j})^{2}\geq \sum\limits_{v_{i}\in V(T_{3})\setminus \{w_{3}\}}\sum\limits_{v_{j}\in V(T_{2})}(x_{v_i}-x_{v_j})^{2}$,
then let $T'_{0}=T_{0}-\sum\limits_{v_{j}\in N_{T_{3}}(w_3)}w_{3}v_{j}+\sum\limits_{v_{j}\in N_{T_{3}}(w_3)}w_{2}v_{j}$,
and otherwise $T'_{0}=T_{0}-\sum\limits_{v_{j}\in N_{T_{3}}(w_3)}w_{3}v_{j}+\sum\limits_{v_{j}\in N_{T_{3}}(w_3)}w_{1}v_{j}$.
Thus, $T'_{0}\in\mathfrak{T}(n,k)$, and by Lemma 3.1, $\rho_{L}(T'_{0})>\rho_{L}(T_{0})$.
This contradicts the maximality of $T_{0}$, and so $n^3_{T_{0}}=2$.\vskip 2mm

This shows that $T_{0}$ consists of a path with its two end vertices $w_1$ and $w_2$ attached some pendant paths respectively.

{\bf Claim 2.} All pendant paths in $T_{0}$ are of length one.

Let $P=w_{1}v_{1}\ldots v_{t}w_{2}$ be the internal path in $T_{0}$.
Suppose on the contrary that $T_{2}=w_{1}u_{1} \ldots u_{s}$ is a pendant path of length $s\geq 2$.
Let $T_{0}=T_{1}\cup T_{2}\cup T_{3}$, where $V(T_{i})\cap V(T_{j})=w_1$ for $1\leq i\neq j\leq 3$, $N_{T_0}(w_{1})\backslash\{u_{1},v_{1}\} \subset V(T_{1})$ and $w_{2}\in V(T_3)$, as in Fig. 2.

If $\sum\limits_{v_{i}\in V(T_{1})}\sum\limits_{v_{j}\in V(T_{3})}(x_{v_i}-x_{v_j})^{2}\geq \sum\limits_{v_{i}\in V(T_{1})}\sum\limits_{v_{j}\in V(T_{2})}(x_{v_i}-x_{v_j})^{2}$,
then let $T'_{0}=T_{0}-\sum\limits_{v_{j}\in V(T_{1})}w_{1}v_{j}+\sum\limits_{v_{j}\in V(T_1)}u_{s-1}v_j$, and otherwise $T'_{0}=T_{0}-\sum\limits_{v_{j}\in V(T_{1})}w_{1}v_{j}+\sum\limits_{v_{j}\in V(T_1)}w_{2}v_j$.
Then $T'_{0}\in\mathfrak{T}(n,k)$ and by Lemma 3.1, $\rho_{Q}(T'_{0})>\rho_{Q}(T_{0})$, which contradicts the maximality of $T_{0}$.

The Claims 1 and 2 show that $ T_{0}\in \mathfrak{D}(n,k)$. $\Box$\vskip 2mm

{\bf Lemma 3.4} \cite{H. Minc}. {\it If $u$ and $v$ are two non-adjacent vertices of graph $G$, then $\rho_{Q} (G+uv)> \rho_{Q}(G)$.}\vskip 2mm

{\bf Theorem 3.5.} {\it The graph with the maximum distance signless Laplacian spectral radius in $\Re(n,k)$ lies to $\mathfrak{D}(n,k)$.}

{\bf Proof.} Suppose that graph $G\in \Re(n,k)$, and that $H$ is a spanning tree of $G$.
Then by Lemma 3.4, $\rho_{Q}(H)>\rho_{Q}(G)$.
Let $\mathfrak{K}(H)$ denote the number of pendant vertices in the tree $H$.
Then $\mathfrak{K}(H)\geq k$.
Theorem 3.3 has told us that $\rho_{Q}(T(n,\mathfrak{K}(H); n_1,n_2))\geq \rho_{Q}(H)$, with equality if and only if $H\cong T(n,\mathfrak{K}(H); n_1,n_2)$,
where $\mathfrak{K}(H)=n_1+n_2$. If $\mathfrak{K}(H)=k$ then the proof is finished,
and otherwise we can construct a tree $T\in \mathfrak{T}(n,k)$ obtained from $T(n,\mathfrak{K}(H); n_1,n_2)$
by deleting $\mathfrak{K}(H)-k$ pendant edges and then attaching them with the remaining pendant edges.
By Lemma 3.2, $\rho_{Q}(T)> \rho_{Q}(T(n,\mathfrak{K}(H); n_1,n_2))$.
By Theorem 3.3, there exist some $T(n,k; t_1,t_2)\in \mathfrak{D}(n,k)$ such that $\rho_{Q}(T(n,k; t_1,t_2))>\rho_{Q}(T)$.
Thus, $\rho_{Q}(T(n,k; t_1,t_2))>\rho_{Q}(G)$. $\Box$\vskip 2mm

As one consequence of the above argument, we directly have\vskip 2mm

{\bf Theorem 3.6.} {\it If $G\in \Re(n,3)$ attains the maximum distance signless Laplacian spectral radius then $G\cong T(n,3; 1,2)$.}


\begin{thebibliography}{10}
\bibitem{M. Aouchiche} M. Aouchiche, P. Hansen, Two Laplacians for the distance
matrix of a graph, Linear Algebra Appl. 439 (2013) 21-33.
\bibitem{M. Aouchiche-1} M. Aouchiche, P. Hansen, Some properties of the distance Laplacian eigenvalues of a graph, Czechoslovak Mathematical Journal 64 (2014) 751-761.
\bibitem{M. Aouchiche-2} M. Aouchiche, P. Hansen, The distance spectra of a graphs: A survey, Linear Algebra Appl. 458 (2014) 301-386.
\bibitem{S. Bose} S. Bose, M. Nath, S. Paul, On the maximal distance spectral radius of graphs without a pendant vertex, Linear Algebra Appl. 438 (2013) 4260-4278.
\bibitem{A. Ilic}  A. Ili\'c, Distance spetral radius of trees with given matching number,
Discrete Appl. Math. 158 (2010) 1799-1806.
\bibitem{H. Liu} H. Liu, M. Lu, Bounds On the distance signless Laplacian spectral radius in terms of clique number, Linear and Multilinaer Algebra 63 (2015) 1750-1759.
\bibitem{H. Lin-2} H. Liu, B. Zhou, On the distance Laplacian spectral radius of graphs, Linear Algebra Appl. 475 (2015) 265-275.
\bibitem{H. Minc} H. Minc, Nonnegative matrices. New York; 1988.
\bibitem{M. Nath} M. Nath, S. Paul, On the distance Laplacian spectra of graphs, Linear Algebra Appl.
460 (2014) 97-110.
\bibitem{W. Ning} W. Ning, L. Ouyang, M. Lu, Distance spectral radius of trees with fixed number of pendant vertices, Linear Algebra Appl. 439 (2013) 2240-2249.
\bibitem{A. Niu} A. Niu, D. Fan and G. Wang, On the distance Laplacian spectra of bipartite graphs, Discrete Appl. Math. 186 (2015) 207-213.
\bibitem{D. Stevanovic} D. Stevanovi\'c, A. Ili\'c, Spectral properties of distance matrix of graphs,
in: I. Gutman, B. Furtula (Eds), Distance in Molecular graphs Theory, in: Math. Chem. Monogr., Vol. 12, University of Kragujevac, 2010, pp. 139-176.
\bibitem{R. Xing} R. Xing, B. Zhou, On the distance and distance signless Laplacian
spectral radii of bicyclic graphs, Linear Algebra Appl. 439 (2013) 3955-3963.
\bibitem{R. Xing-2} R. Xing, B. Zhou, J. Li, On the distance signless Laplacian spectral radius of graphs,
Linear and Multilinear Algebra 62 (2014) 1377-1387.
\bibitem{G. Yu} G. Yu, H. Jia, H. Zhang, J. Shu, Some graft transformations and its applications on
the distance spectral radius of a graph, Appl. Math. Letters 25 (2012) 315-319.
\bibitem{G. Yu-2} G. Yu, Y. Wu, Y. Zhang, J. Shu, Some graft transformations and its application
on a distance spectrum, Discrete Math. 311 (2011) 2117-2123.
\end{thebibliography}
\end{document}